\documentclass[10pt]{amsart}
\usepackage{amssymb}
\usepackage[headings]{fullpage}
\usepackage{amsfonts}
\usepackage{paralist}
\usepackage{color}

\usepackage[colorlinks,pagebackref=true,pdftex]{hyperref}

\makeatletter
\@addtoreset{equation}{section}
\def\theequation{\thesection.\@arabic \c@equation}

\def\theenumi{\@roman\c@enumi}
\def\@citecolor{blue}
\def\@linkcolor{blue}
\def\@urlcolor{blue}

\makeatother

\newtheorem{claim*}{Claim}

\theoremstyle{definition}
\newtheorem{remark}[equation]{Remark}

\newtheorem{eg}[equation]{Example}

\newtheorem{definition}[equation]{Definition}

\newtheorem{notn}[equation]{Notation}

\newtheorem{innercustomgeneric}{\customgenericname}
\providecommand{\customgenericname}{}
\newcommand{\newcustomtheorem}[2]{%
  \newenvironment{#1}[1]
  {%
   \renewcommand\customgenericname{#2}%
   \renewcommand\theinnercustomgeneric{##1}%
   \innercustomgeneric
  }
  {\endinnercustomgeneric}
}

\newcustomtheorem{customthm}{Theorem}
\newcustomtheorem{customlemma}{Lemma}

\newsavebox{\upperboundtheorem}

\title{The Pentagonal Inequality}

\author[R.~Barbara]{Roy Barbara}

\begin{document}

\maketitle

\section{ \textbf{Introduction}}

We recall the following result, conjectured by L. Fejes T\'{o}th (see [\ref{ref1}]) and proved by Lenhard H-C (see [\ref{ref2}]).

\begin{customthm}{0}\label{thm0}(T\'{o}th's inequality) Let $n\geq3$. Let $x_1, x_2, ..., x_n$ and $\alpha_1, \alpha_2, ..., \alpha_n$ be positive real numbers with $\sum\limits_{i=1}^n\alpha_i=\pi.$ Then, with $x_{n+1}=x_1$, we have the following inequality:  
$$\sum\limits_{i=1}^nx_ix_{i+1}cos\alpha_i\hspace{0.1cm}\leq\hspace{0.1cm} cos\frac{\pi}{n}.\sum\limits_{i=1}^nx_i^2$$ Our interest lies in a sharp majoration of a sum as $\sum\limits_{i=1}^na_icos\alpha_i,$ where $a_i, \alpha_i>0$ and $\sum\limits_{i=1}^n\alpha_i=\pi.$ Results for the cases $n=3, 4$ can be found in [\ref{ref3}]. In this article, we focus our attention on the case $n=5$. In section \ref{sec4}, we highlight the sharpness of the result. In section \ref{sec5}, we quickly study the case $n=7$. The case $n=6$ remains open. 
\end{customthm}

\section{ \textbf{The Results}}\label{sec2}

We present the Pentagonal Inequality in two forms: the strong form (theorem \ref{thm1}) and the normal form (theorem \ref{thm2}). Though theorem \ref{thm1} is more precise, theorem \ref{thm2} remains sometimes more practical.

\begin{customthm}{1}\label{thm1}
Let $a_1, a_2, ..., a_5$ and $\alpha_1, \alpha_2, ..., \alpha_5$ be positive real numbers with $\sum\limits_{i=1}^5\alpha_i=\pi$ and such that $a_1\leq a_2\leq a_3\leq a_4\leq a_5.$ \hspace{0.1cm}  Then, we have the sharp inequality: $$a_1cos\alpha_1+a_2cos\alpha_2+a_3cos\alpha_3+a_4cos\alpha_4+a_5cos\alpha_5$$ $$\leq$$ 
$$\big{(}\frac{1+\sqrt5}{4})\hspace{0.05cm}\frac{1}{a_1a_2a_3a_4a_5}\hspace{0.05cm}(a_1^2a_5^2a_2^2+a_5^2a_2^2a_3^2+a_2^2a_3^2a_4^2+a_3^2a_4^2a_1^2+a_4^2a_1^2a_5^2)$$
\end{customthm}

\begin{customthm}{2}\label{thm2}
Let $a_1, a_2, ..., a_5$ and $\alpha_1, \alpha_2, ..., \alpha_5$ be positive real numbers with $\sum\limits_{i=1}^5\alpha_i=\pi$. Then, we have: 
$$a_1cos\alpha_1+a_2cos\alpha_2+a_3cos\alpha_3+a_4cos\alpha_4+a_5cos\alpha_5$$
$$\leq$$ 
$$(\frac{1+\sqrt5}{4})\hspace{0.05cm}\frac{1}{a_1a_2a_3a_4a_5}\hspace{0.05cm}(a_1^2a_2^2a_3^2+a_2^2a_3^2a_4^2+a_3^2a_4^2a_5^2+a_4^2a_5^2a_1^2+a_5^2a_1^2a_2^2)$$
\end{customthm}
\vspace{1cm}

\noindent We briefly indicate that the case of equality in either theorem  \ref{thm1} or \ref{thm2} requires the condition: 
$$a_1sin\alpha_1=a_2sin\alpha_2=a_3sin\alpha_3=a_4sin\alpha_4=a_5sin\alpha_5.$$

 \section{ \textbf{Proofs of Theorems 1 and 2}}\label{sec4}

Let $a_1, a_2, ..., a_5$ be positive real numbers. A circular permutation of the $a_i$ can be represented by $\sigma=(x_1, x_2, ..., x_5)$, with $x_1=a_1$, and where $x_2, x_3, x_4, x_5$ are $a_2, a_3, a_4, a_5$ in some order. There are $4!=24$ such permutations. 

\noindent {\itshape{In all what follows}}, we use the following notation:

\noindent \textbf{Notation}: Let $\sigma=(x_1, x_2, ..., x_5)$ be a circular permutation of $a_1, a_2, ..., a_5$ (with $x_1=a_1$). We define $\phi(\sigma)$, also denoted by $\phi(x_1, x_2, ..., x_5)$ as follows: $$\phi(x_1, x_2, x_3, x_4, x_5)=x_1x_2x_3+x_2x_3x_4+x_3x_4x_5+x_4x_5x_1+x_5x_1x_2.$$

\begin{customlemma}{1}\label{lemma1}
Let $a_1\leq a_2\leq ...\leq a_5$ be positive real numbers. Then, there is a circular permutation of the $a_i$,  say $\sigma_0$, \hspace{0.05cm}such that $\phi(\sigma_0)$ minimizes all the $\phi(\sigma)$. Namely, $\sigma_0$ is the cycle\hspace{0.1cm} $\sigma_0=(a_1, a_5, a_2, a_3, a_4).$

\noindent Equivalently, one might take $\sigma_1=(a_1, a_4, a_3, a_2, a_5)$ \hspace{0.1cm} since $\phi(\sigma_1)=\phi(\sigma_0).$
\end{customlemma}

\begin{proof}
After some amount of algebra, we find the following identities (and the result follows):
\begin{align*}
\phi(a_1, a_2, a_3, a_4, a_5)-\phi(\sigma_0) & = \phi(a_1, a_5, a_4, a_3, a_2)-\phi(\sigma_0)\\ &  =a_3(a_4-a_2)(a_5-a_1)\geq0
\end{align*}
\begin{align*}
\phi(a_1, a_2, a_3, a_5, a_4)-\phi(\sigma_0) &= \phi(a_1, a_4, a_5, a_3, a_2)-\phi(\sigma_0)\\ &=a_1(a_3-a_2)(a_5-a_4)+a_3(a_4-a_1)(a_5-a_2)\geq 0
\end{align*}
\begin{align*}
\phi(a_1, a_2, a_4, a_3, a_5)-\phi(\sigma_0)&= \phi(a_1, a_5, a_3, a_4, a_2)-\phi(\sigma_0)\\ &=a_1(a_3-a_2)(a_5-a_4)+a_5(a_3-a_1)(a_4-a_2)\geq 0
\end{align*}
\begin{align*}
& \phi(a_1, a_2, a_4, a_5, a_3)-\phi(\sigma_0)= \phi(a_1, a_3, a_5, a_4, a_2)-\phi(\sigma_0)\\ &=a_1(a_3-a_2)(a_5-a_4)+a_3(a_4-a_1)(a_5-a_2)+a_5(a_2-a_1)(a_4-a_3)\geq 0
\end{align*}
\begin{align*}
\phi(a_1, a_2, a_5, a_3, a_4)-\phi(\sigma_0)&= \phi(a_1, a_4, a_3, a_5, a_2)-\phi(\sigma_0)\\ &=a_4(a_3-a_1)(a_5-a_2)\geq 0
\end{align*}
\begin{align*}
\phi(a_1, a_2, a_5, a_4, a_3)-\phi(\sigma_0)&= \phi(a_1, a_3, a_4, a_5, a_2)-\phi(\sigma_0)\\ &=a_3(a_4-a_1)(a_5-a_2)+a_5(a_2-a_1)(a_4-a_3)\geq 0
\end{align*}
\begin{align*}
\phi(a_1, a_3, a_2, a_4, a_5)-\phi(\sigma_0)&= \phi(a_1, a_5, a_4, a_2, a_3)-\phi(\sigma_0)\\ &=a_1(a_3-a_2)(a_5-a_4)+a_2(a_4-a_3)(a_5-a_1)\geq 0
\end{align*}
\begin{align*}
\phi(a_1, a_3, a_2, a_5, a_4)-\phi(\sigma_0)&= \phi(a_1, a_4, a_5, a_2, a_3)-\phi(\sigma_0)\\ &=a_2(a_4-a_1)(a_5-a_3)\geq 0
\end{align*}
\begin{align*}
\phi(a_1, a_3, a_4, a_2, a_5)-\phi(\sigma_0)&= \phi(a_1, a_5, a_2, a_4, a_3)-\phi(\sigma_0)\\ &=a_5(a_2-a_1)(a_4-a_3)\geq 0
\end{align*}
\begin{align*}
\phi(a_1, a_3, a_5, a_2, a_4)-\phi(\sigma_0)&= \phi(a_1, a_4, a_2, a_5, a_3)-\phi(\sigma_0)\\ &=a_1(a_3-a_2)(a_5-a_4)+a_4(a_2-a_1)(a_5-a_3)\geq 0
\end{align*}
\begin{align*}
\phi(a_1, a_4, a_2, a_3, a_5)-\phi(\sigma_0)= \phi(a_1, a_5, a_3, a_2, a_4)-\phi(\sigma_0)=a_1(a_3-a_2)(a_5-a_4)\geq 0\end{align*}
\begin{align*}\phi(a_1, a_4, a_3, a_2, a_5)-\phi(\sigma_0)= \phi(a_1, a_5, a_2, a_3, a_4)-\phi(\sigma_0)=0\end{align*}
\end{proof}

\begin{customlemma}{2}\label{lemma2}
Let $a_1, a_2, ..., a_5$ and $\alpha_1, \alpha_2, ..., \alpha_5$ be positive real numbers with $\sum\limits_{i=1}^5 \alpha_i=\pi$. If $\sigma =(b_1, b_2, ..., b_5)$ is any circular permutation of the $a_i$ (with $b_1=a_1$), then we have: 
\begin{equation}\tag{1}\label{eqn1}
\sum\limits_{i=1}^5 b_icos\alpha_i\hspace{0.1cm}\leq\hspace{0.1cm}\displaystyle{ \big(\frac{1+\sqrt5}{4})\hspace{0.05cm}\frac{1}{b_1b_2b_3b_4b_5}\hspace{0.05cm}\phi(b_1^2b_2^2b_3^2b_4^2b_5^2)}
\end{equation}
\end{customlemma}

\begin{proof}
Set $\beta_1=\alpha_1,\hspace{0.1cm} \beta_2=\alpha_4,\hspace{0.1cm} \beta_3=\alpha_2,\hspace{0.1cm} \beta_4=\alpha_5,$ and $\beta_5=\alpha_3$ (so $\sum\limits_{i=1}^5\beta_i=\pi).$

\noindent Set $P=b_1b_2b_3b_4b_5=a_1a_2a_3a_4a_5,$\hspace{0.05cm} and set\vspace{0.2cm}\\ \noindent $\displaystyle{x_1=\sqrt{\frac{b_1b_2b_3}{b_4b_5}},\hspace{0.3cm} x_2=\sqrt{\frac{b_1b_4b_5}{b_2b_3}},\hspace{0.3cm}  x_3=\sqrt{\frac{b_2b_3b_4}{b_1b_5}},\hspace{0.3cm} x_4=\sqrt{\frac{b_1b_2b_5}{b_3b_4}}},\hspace{0.3cm}  x_5=\sqrt{\frac{b_3b_4b_5}{b_1b_2}}.$

\vspace{0.2cm}
\noindent We have $x_1^2=P^{-1}(b_1^2b_2^2b_3^2),\hspace{0.1cm} x_3^2=P^{-1}(b_2^2b_3^2b_4^2), \hspace{0.1cm}x_5^2=P^{-1}(b_3^2b_4^2b_5^2),$ \hspace{0.1cm}$x_2^2=P^{-1}(b_4^2b_5^2b_1^2),$ and $x_4^2=P^{-1}(b_5^2b_1^2b_2^2)$. Summing, we obtain: \begin{equation}\tag{2}\label{eqn2} \sum\limits_{i=1}^5x_i^2 =P^{-1}\hspace{0.05cm}\phi(b_1^2, b_2^2, b_3^2, b_4^2, b_5^2)\end{equation}
On the other hand we have: $x_1x_2cos\beta_1=b_1cos\alpha_1,\hspace{0.1cm} x_2x_3cos\beta_2=b_4cos\alpha_4,\hspace{0.1cm}$ $x_3x_4cos\beta_3=b_2cos\alpha_2,\hspace{0.1cm}$ $x_4x_5cos\beta_4=b_5cos\alpha_5,$ \hspace{0.1cm}and\hspace{0.1cm} $x_5x_1cos\beta_5=b_3cos\alpha_3.$ 

\noindent Hence, with $x_6=x_1$, we get: \begin{equation}\tag{3}\label{eqn3}\sum\limits_{i=1}^5 x_ix_{i+1}cos\beta_i=\sum\limits_{i=1}^5 b_icos\alpha_i\end{equation}   
\noindent By theorem \ref{thm0}, we have $\sum\limits_{i=1}^5 x_ix_{i+1}cos\beta_i\leq cos\frac{\pi}{5}.\sum\limits_{i=1}^5 x_i^2$. From this, (\ref{eqn2}), (\ref{eqn3}), and $\displaystyle{cos\frac{\pi}{5}=\frac{1+\sqrt5}{4}}$, we obtain (\ref{eqn1}). 

\end{proof}
\noindent {\it{Proofs of theorems \ref{thm1} and \ref{thm2}}}: To get theorem \ref{thm2}, just apply lemma \ref{lemma2} to $\sigma=(a_1, a_2, a_3, a_4, a_5).$ To get theorem \ref{thm1}, apply lemma \ref{lemma2} to the cycle $\sigma_0=(a_1, a_5, a_2, a_3, a_4),$ obtaining, in virtue of lemma \ref{lemma1}, the sharpest form of the Pentagonal Inequality. 

\noindent Note the following: If $a_1< a_2< ...<a_5$ (or more simply if $a_1\leq a_2\leq ...\leq a_5$ and $a_2<a_4)$, then, $\phi(a_1^2, a_5^2, a_2^2, a_3^2, a_4^2)<\phi(a_1^2, a_2^2, a_3^2, a_4^2, a_5^2).$ Hence, {\it{in this case}}, while equality in theorem \ref{thm1} might be reached, equality in theorem \ref{thm2} never arises.

\vspace{1cm}
\section{ \textbf{Sharpness of the Pentagonal Inequality}}\label{sec4}

To see why the Pentagonal Inequality is {\it{sharp}}, we prove that T\'{o}th's inequality for $n=5$ can be {\it{derived}} from the Pentagonal Inequality (the normal form suffices): Indeed, let $x_1, x_2, ..., x_5$ and $\alpha_1, \alpha_2, ..., \alpha_5$ be positive real numbers with $\sum\limits_{i=1}^5\alpha_i=\pi$. Given by hypothesis that theorem \ref{thm2} holds, we prove (with $x_6=x_1$) that \begin{equation}\tag{4}\label{eqn4}\sum\limits_{i=1}^5 x_ix_{i+1}cos\alpha_i\hspace{0.1cm}\leq\hspace{0.1cm} cos\frac{\pi}{5}.\sum\limits_{i=1}^5 x_i^2 \end{equation}
Set $\beta_1=\alpha_1,\hspace{0.1cm} \beta_2=\alpha_3,\hspace{0.1cm} \beta_3=\alpha_5,\hspace{0.1cm} \beta_4=\alpha_2,$ \hspace{0.07cm}and \hspace{0.07cm}$\beta_5=\alpha_4$ (so $\sum\limits_{i=1}^5 \beta_i =\pi$). Set also $a_1=x_1x_2,\hspace{0.1cm} a_2=x_3x_4,\hspace{0.1cm}$ $a_3=x_5x_1, \hspace{0.1cm}a_4=x_2x_3,$ \hspace{0.07cm}and\hspace{0.07cm} $a_5=x_4x_5.$

\noindent By theorem \ref{thm2}, we have: \begin{equation}\tag{5}\label{eqn5} \sum\limits_{i=1}^5 a_icos\beta_i\hspace{0.1cm}\leq\hspace{0.1cm}\displaystyle{ \big(\frac{1+\sqrt5}{4})\hspace{0.05cm}\frac{1}{a_1a_2a_3a_4a_5}\hspace{0.05cm}\phi(a_1^2a_2^2a_3^2a_4^2a_5^2)}\end{equation}

\noindent We have: $a_1cos\beta_1=x_1x_2cos\alpha_1,\hspace{0.1cm} a_2cos\beta_2=x_3x_4cos\alpha_3,$ $a_3cos\beta_3=x_5x_1cos\alpha_5,\hspace{0.1cm} a_4cos\beta_4=x_2x_3cos\alpha_2,$\hspace{0.1cm} and $a_5cos\beta_5=x_4x_5cos\alpha_4.$\hspace{0.1cm} {\it{Hence}}, \begin{equation}\tag{6}\label{eqn6} \sum\limits_{i=1}^5 a_icos\beta_i=\sum\limits_{i=1}^5x_ix_{i+1}cos\alpha_i\end{equation}

\noindent Set \hspace{0.1cm}$P=a_1a_2a_3a_4a_5.$ \hspace{0.15cm}We easily find that \hspace{0.15cm}$P^{-1}\vspace{0.07cm}(a_1^2a_2^2a_3^2)=x_1^2,\hspace{0.15cm} P^{-1}(a_2^2a_3^2a_4^2)=x_3^2, \hspace{0.15cm}P^{-1}(a_3^2a_4^2a_5^2)=x_5^2,$\\ \noindent$P^{-1}(a_4^2a_5^2a_1^2)=x_2^2,$ and $P^{-1}(a_5^2a_1^2a_2^2)=x_4^2.$ 

\noindent Summing yields $\sum\limits_{i=1}^5 x_i^2=P^{-1}\phi(a_1^2, a_2^2, a_3^2, a_4^2, a_5^2).$ Hence the right-side in (5) is $cos\frac{\pi}{5}.\sum\limits_{i=1}^5x_i^2.$ From this, (\ref{eqn5}), and (\ref{eqn6}), we get (\ref{eqn4}).

\vspace{1cm}
\section{ \textbf{The Heptagonal Inequality}}\label{sec5}

The case of an odd integer $n\geq 7$ can be treated similarly. For $n=7$, we obtain: 

\begin{customthm}{3}\label{thm3}(the normal form) Let $a_1, a_2, ..., a_7$ and $\alpha_1, \alpha_2, ..., \alpha_7$ be positive real numbers with $\sum\limits_{i=1}^7\alpha_i=\pi.$ Then, we have: $$\sum\limits_{i=1}^7a_icos\alpha_i\hspace{0.1cm}\leq\hspace{0.1cm} cos\frac{\pi}{7}.\frac{1}{a_1a_2...a_7}\hspace{0.05cm} \psi(a_1^2, a_2^2, a_3^2, a_4^2, a_5^2, a_6^2, a_7^2)$$
where \begin{align*} \psi(x_1, x_2, ..., x_7)=x_1x_2x_3x_4& +x_2x_3x_4x_5+x_3x_4x_5x_6+x_4x_5x_6x_7\\ &+x_5x_6x_7x_1+x_6x_7x_1x_2+x_7x_1x_2x_3.\end{align*}
\end{customthm}

\begin{proof}
Set $\beta_1=\alpha_1,\hspace{0.1cm} \beta_2=\alpha_5, \hspace{0.1cm}\beta_3=\alpha_2, \hspace{0.1cm}\beta_4=\alpha_6, \hspace{0.1cm}\beta_5=\alpha_3, \hspace{0.1cm}\beta_6=\alpha_7, \hspace{0.1cm}\beta_7=\alpha_4,$ \hspace{0.05cm}and set
\vspace{0.3cm}\\
\noindent $\displaystyle{x_1=\sqrt{\frac{a_1a_2a_3a_4}{a_5a_6a_7}}, \hspace{0.15cm}x_2=\sqrt{\frac{a_1a_5a_6a_7}{a_2a_3a_4}}, \hspace{0.15cm}x_3=\sqrt{\frac{a_2a_3a_4a_5}{a_1a_6a_7}}, \hspace{0.15cm}x_4=\sqrt{\frac{a_1a_2a_6a_7}{a_3a_4a_5}}}, $
$\displaystyle{\hspace{0.15cm}x_5=\sqrt{\frac{a_3a_4a_5a_6}{a_1a_2a_7}}},$\\

\vspace{0.2cm}\noindent \vspace{0.2cm}$\vspace{0.2cm}\displaystyle{x_6=\sqrt{\frac{a_1a_2a_3a_7}{a_4a_5a_6}}, \hspace{0.1cm} \text{and }  \hspace{0.1cm}x_7=\sqrt{\frac{a_4a_5a_6a_7}{a_1a_2a_3}}}.$

\noindent Then, apply theorem \ref{thm0} with these values of $x_1, x_2, ..., x_7$ and $\beta_1, \beta_2, ..., \beta_7.$
\end{proof}

\providecommand{\bysame}{\leavevmode\hbox to3em{\hrulefill}\thinspace}

\providecommand{\MRhref}[2]{%
  \href{http://www.ams.org/mathscinet-getitem?mr=#1}{#2}
}
\providecommand{\href}[2]{#2}

\end{document}